\documentclass{ifacconf}

\usepackage[square,sort,comma,numbers]{natbib}
\usepackage{graphics}
\usepackage{graphicx} 

\usepackage{epstopdf}  
\usepackage{epsfig} 
\usepackage{times} 
\usepackage{amsmath} 
\usepackage{amssymb}  
\usepackage{float} 

\usepackage{color}
\usepackage{caption}

\usepackage{mathtools}

\makeatletter
\newcommand{\figcaption}{\def\@captype{figure}\caption}
\newcommand{\tabcaption}{\def\@captype{table}\caption}
\makeatother
\begin{document}
\begin{frontmatter}

\title{
\hskip -.109in  On Drawdown-Modulated Feedback \\ Control in Stock Trading \hskip -.109in
} 


\author{\large Chung-Han Hsieh and B. Ross Barmish}

\vskip -.1in
\address{ECE Department\\ University of Wisconsin\\
	Madison, WI 53706\\
	e-mail: hsieh23@wisc.edu, barmish@engr.wisc.edu}

\begin{abstract}  
Control of drawdown, that is, the control of the drops in wealth over time from peaks to subsequent lows, is of great concern from a risk management perspective. 
With this motivation in mind, the focal point of this paper is to address the drawdown issue in a stock trading context. Although our analysis can be carried out without reference to control theory, to make the work accessible to this community, we use the language of feedback systems. The takeoff point for the results to follow, which we call the {\it Drawdown Modulation Lemma}, characterizes {\it any} investment which guarantees that the percentage  drawdown is no greater than a prespecified level with probability one. With the aid of this lemma, we introduce a new scheme which we call the {\it drawdown-modulated feedback control}.  To illustrate the power of the theory,  we consider a drawdown-constrained version of the well-known Kelly Optimization Problem which involves maximizing the expected logarithmic growth of the trader's account value. As the drawdown parameter $d_{\max}$ in our new formulation tends to one, we recover existing results as a special case.
This new theory leads to an optimal investment strategy whose application is illustrated via an example with historical stock-price data.
\end{abstract}


\end{frontmatter}
\vspace{2mm}
\section{Introduction}
\label{Section1:introduction}
\vspace{-3mm}
Control of drawdown, that is, the control of the drops in wealth over time from peaks to subsequent lows, is of great concern from a risk management perspective. Suffice it to say, the issue of drawdown control has received a considerable attention in the finance literature. More specifically, to properly control drawdown in a stock trading context, one standard approach is to incorporate this consideration into a constrained optimization framework which typically includes other performance criteria;~e.g.,~see~\cite{Grossman_Zhou_1993} and \cite{Klass_Nowicki_2005} which focus on a single-stock scenario. 
 There are also some papers dealing with modifications and extensions of these results for the single-stock case to address an entire portfolio under various stochastic modeling assumptions;~e.g.,~see~\mbox{\cite{Zhou_Shang_2015}-\cite{Cherney_2015}.} 
Finally, the literature also includes various methodologies to address different types of drawdown.
 For example, absolute drawdown is studied in~\cite{Hayes} and other drawdown-based metrics are considered in~\cite{Barmish_Hsieh_2015},~\cite{Rockafellar_2006}~and~\cite{Goldberg_Mahmoud_2016}.

\vspace{5mm}
In this paper, we focus on the notion of {\it percentage drawdown} whose technical definition is given in the next section. We provide a new result which enables a trader to guarantee that a prescribed maximum level $d_{\max}$  for this quantity will not be exceeded. 
To provide further context for this paper, we mention a sampling of some other papers in the existing literature using different {\it risk} measures rather than drawdown. Examples of such measures include Value at Risk~(VaR), Conditional Value at~Risk, Expected Shortfall and the celebrated mean-variance criterion;~e.g.,~see~\mbox{\cite{Jorion_2006}-\cite{Luenberger_2011}} and \cite{Markowitz_1959}.
In addition to the analysis of risk, some papers in the literature consider portfolio optimization involving a maximization of expected logarithmic growth. This is the so-called Kelly Optimization Problem which will be used to demonstrate our theory; e.g.,~\mbox{see~\cite{Maclean_1992}-\cite{MacLean_Thorp_Ziemba_2011}. }
Related to this, the literature also includes a well-known method called the {\it Fractional Kelly Strategy}. This is aimed at scaling down the size of investment for the purpose of mitigating  the risk;~e.g., see~\cite{Thorp_2006} and~\cite{Maclean_Thorp_Ziemba_2010}.   

\vspace{5mm}
The takeoff point for this paper, which we call the {\it Drawdown Modulation Lemma}, characterizes investments which guarantee that the  {\it percentage drawdown} is no greater than a prespecified level~$d_{\max} \in (0,1)$ with probability one. To make our exposition appealing to this community, this investment scheme which we derive from the lemma, is expressed in a classical feedback control setting. We call it {\it drawdown-modulated feedback control}.
As the drawdown parameter $d_{\max}$, we recover existing results as a special case. To further illustrate its use, as previously mentioned, we consider a drawdown-constrained version of the Kelly Optimization Problem which involves maximizing the expected logarithmic growth.  Then, a numerical example with historical data is carried out and a further generalization of the lemma for portfolio case is discussed in this paper.

\vspace{8mm}
 \section{Drawdown Definitions}
 \vspace{-3mm} 
For $k=0,1,2,...,N$, we let $V(k)$ denote the account value at stage~$k$. Then as $k$ evolves, the {\it percentage drawdown}~(to-date) is defined as 
 \[
 d(k) \doteq \frac{V_{\max}(k) - V(k)}{V_{\max}(k)}
 \] where 
 $$
 V_{\max}(k) \doteq \max_{0 \le i \le k}V(i).
 $$
This leads to {\it overall percentage drawdown}
 $$ 
 d_{\max}^* \doteq \max_{0 \leq k \leq N} d(k).
 $$
Note that the percentage drawdown satisfies~$0 \leq d(k) \leq 1.$ Although not considered here, there is another well-known measure, called the {\it maximum absolute drawdown}, which we denote by~$D_{\max}^*$ and is given by
\[
D_{\max}^* \doteq \max_{0\leq k\leq N} V_{\max}(k) - V(k).
\]
The reader is referred to~\cite{Ismail_2004} and~\cite{Malekpour_Barmish_2013} for work on this topic.

\vspace{5mm}
To further elaborate on these two notions of drawdown, we consider an example with a hypothetical account value~$V(k)$ shown in Figure~\ref{fig:Drawdown_Elaboration}.  From the plot, we obtain the overall percentage drawdown, 
 $
 d_{\max}^* \approx 0.8333,
 $
 occurs at stage~$k=7.$ On the other hand, the maximum absolute drawdown,
$
D_{\max}^* = 4.5
$,
 occurs at stage $k=24.$ This toy example shows that the two types of drawdown described above can be quite different. In this paper, we concentrate on the percentage drawdown. This is the version of drawdown which is better suited to deal with different ``scales" for $V(k)$. That is both small and large investors can identify with this concept.
 
\vspace{2mm}
\begin{center}
	\graphicspath{{Figs/}}
	\includegraphics[scale=0.42]{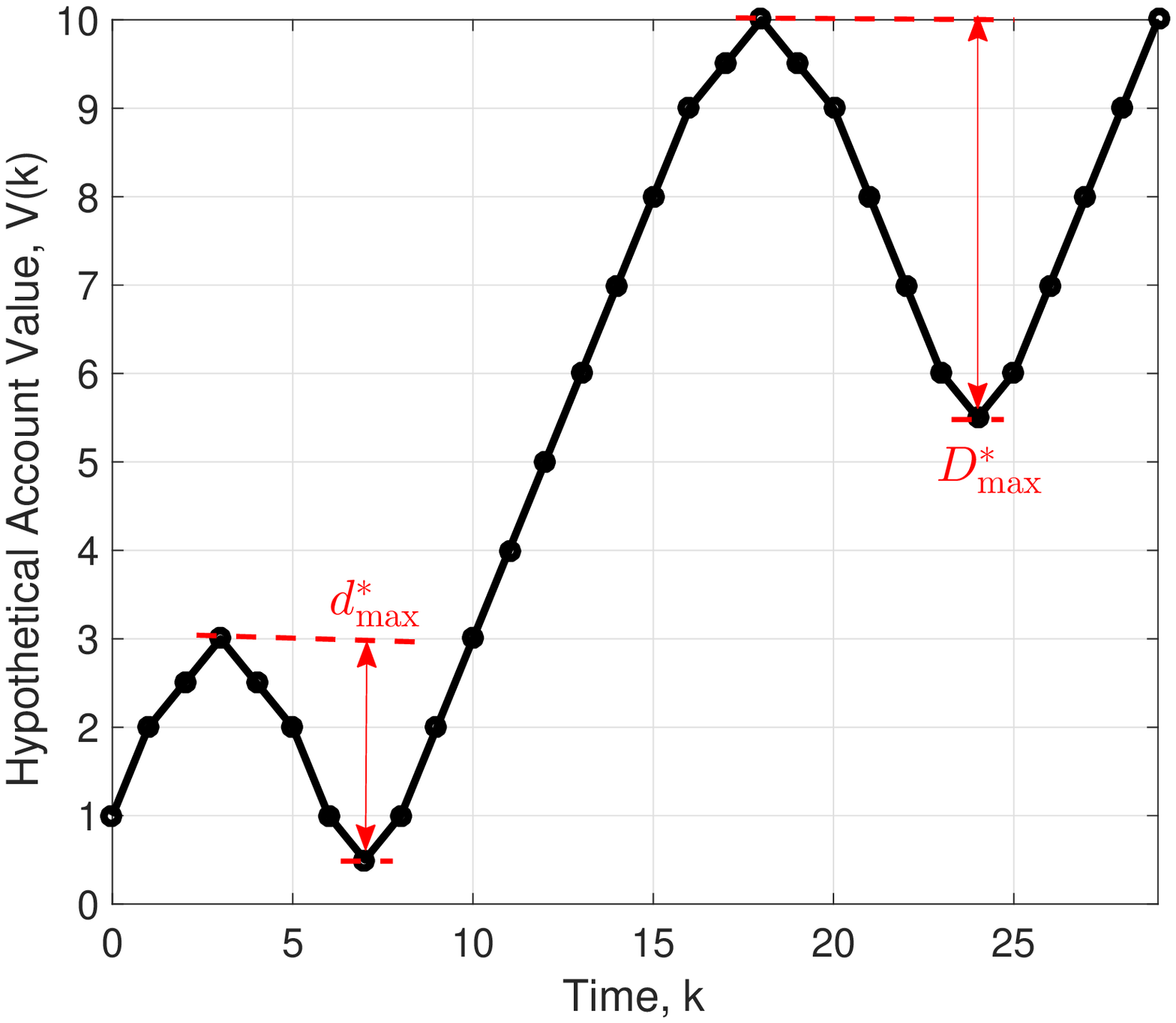}
	\figcaption{Maximum Percentage and Absolute Drawdown}
	\label{fig:Drawdown_Elaboration}
\end{center}

\vspace{8mm}
\section{Feedback Control Point of View} 
\label{Section: Feedback Control Point of View}
\vspace{-3mm}
In the sequel, as previously mentioned, we emphasize the control-theoretic point of view. Our formulation here is consistent with a  growing body of the literature addressing finance problems but originating from the control community; e.g., see~\cite{Barmish_2011}-\cite{Barmish_Primbs_2015}. Although our analysis to follow can be carried out without reference to control theory, to make the work accessible to this audience, we use the language of feedback systems.
Specifically, we view the stock prices~$S(k)$ as exogenous inputs to a feedback system. In this setting, we use a linear feedback controller which modifies the investment~$I(k)$ using a time-varying gain $K(k)$ applied to the account value $V (k)$. That is, for each stage $k$, we consider the {\it controller} having the form
$$
I(k) = K(k) V(k).
$$
Typically, when selecting this feedback gain, we include a constraint which we express as $K(k) \in \mathcal{K}$. For instance, suppose we restrict attention to a so-called {\it cash-financed} investment. Then we impose a constraint which guarantees~$|I(k)| \leq V(k)$. That is, the investment level is limited to the value of one's account. This in turn forces~$|K(k)|\leq 1$. The feedback control configuration which describes this scheme is depicted in Figure~\ref{fig:BlockDiagram} for a single stock. 
Such a configuration can be generalized to deal with a portfolio case of $n$ stocks;~e.g., see~\cite{Barmish_Hsieh_2015}. To link back to finance concepts, the special case of {\it buy-and-hold} is obtained when~$K(k) \equiv 1.$ 

\vspace{5mm}

\begin{center}
	\graphicspath{{Figs/}}
	\includegraphics[scale=0.45]{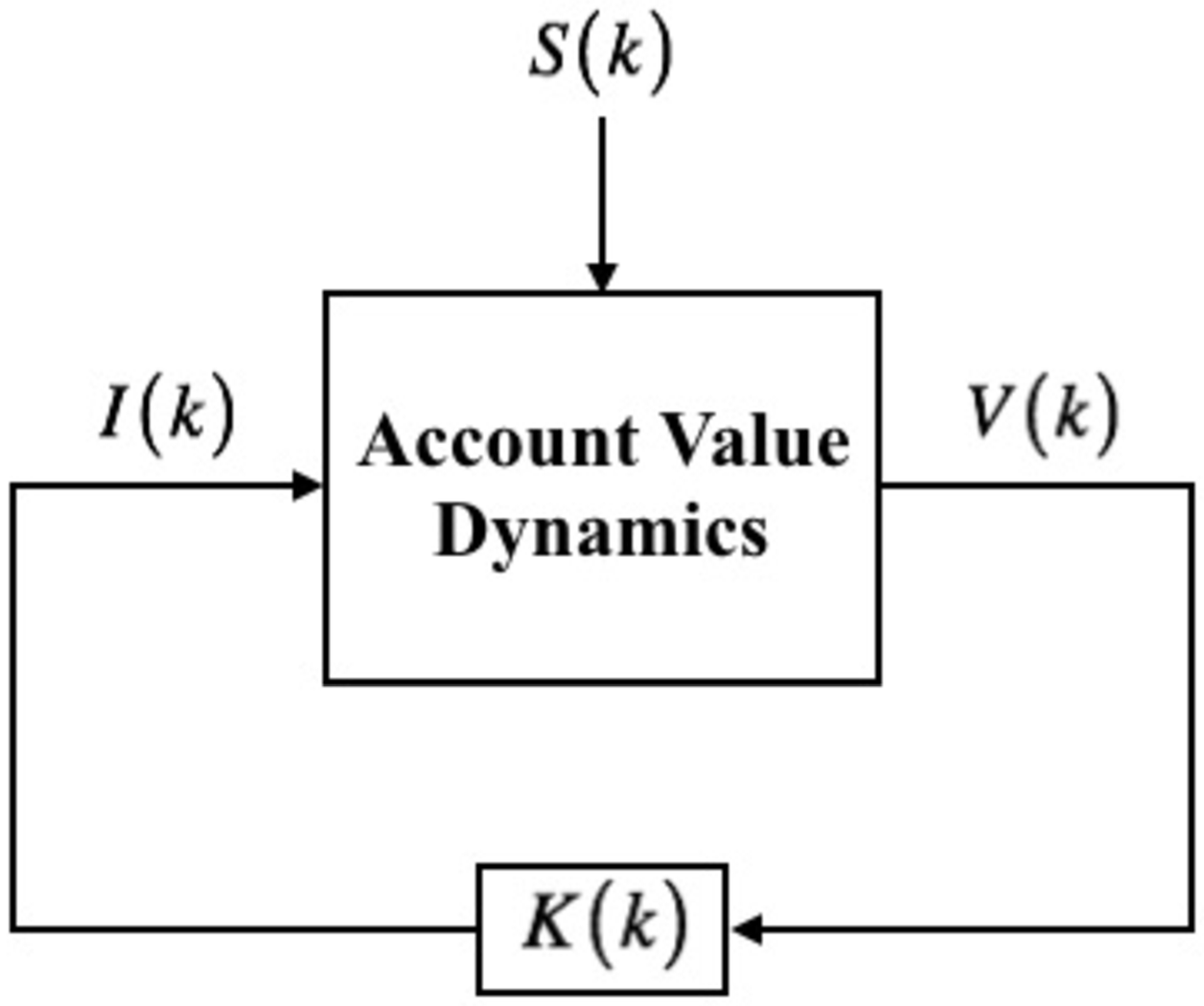}
	\figcaption{Feedback Control Configuration}
	\label{fig:BlockDiagram}
\end{center}


\vspace{8mm}
\section{Stock-Trading Formulation}
\label{Section 2: Problem Formulation}
\vspace{-2mm}
For $k=0,1,2,...,N$, we let~$S(k) >0 $ denote the stock price. The associated returns up to stage $N-1$ are given by
\[
X(k) \doteq \frac{S(k+1) -S(k)}{S(k)}
\]
In the sequel, we assume that the returns are bounded; i.e.,
\[
X_{\min} \leq X(k) \leq X_{\max}
\]
with $X_{\min}$ and $X_{\max}$ being points in the support, denoted by~$\mathcal{X}$, and satisfying
$$
-1 < X_{\min} < 0 < X_{\max}.
$$
In the sections to follow, when necessary, we assume the random variables $X(k)$ are independent and identically distributed~(i.i.d.). 

\vspace{5mm}
{\bf Account Value and Drawdown}:
 Now letting~$I(k)$ be the investment  at stage~$k$ and note that~\mbox{$I(k) < 0$} stands for short selling. Beginning at some initial account value~$V(0)>0$, the evolution to terminal state~$V(N)$ is described sequentially by the recursion
\[
V(k+1) = V(k) + I(k) X(k).
\]
Now, given a {\it maximum acceptable drawdown  level}~$d_{\max}$ satisfying~\mbox{$0 < d_{\max} < 1$}, we focus on conditions on $I(k)$ under which  satisfaction of the constraint 
$$
d(k) \leq d_{\max}
$$
is assured along all sample pathes  with probability one. 

\vspace{5mm}
{\bf Idealized Market}:
In the sequel, we further assume that our stock-trading occurs within ``idealized market." That is, we assume zero transaction costs, zero interest rates and perfect liquidity conditions. For more details about idealized market assumption, the reader is referred to reference~\cite{Barmish_Primbs_2015}.

\vspace{8mm}
\section{The Drawdown Modulation Lemma}
\label{Section 3: The Drawdown Modulation Theorem}
\vspace{-2mm}
In this section, the stepping stone in this paper, which we call it {\it The Drawdown Modulation Lemma}, is given. This lemma provides a necessary and sufficient condition on the investment~$I(k)$ which guarantees that the percentage drawdown is no greater than a given level $d_{\max}$  with probability one. 

\vspace{4mm}
{\bf The Drawdown Modulation Lemma}: {\it An investment function~$I(\cdot)$ guarantees maximum acceptable drawdown level~$d_{\max}$ or less with probability one if and only if for all $k$, the~\mbox{condition}
	\[ 
	- \frac{{{d_{\max }} - d(k)}}{{\left( {1 - d\left( k \right)} \right){X_{\max }}}}V(k) \leq I(k) \leq \frac{{{d_{\max }} - d(k)}}{{\left( {1 - d\left( k \right)} \right)\left| {{X_{\min }}} \right|}}V(k)
	\]
is satisfied along all sample pathes.}

\vspace{5mm}
{\bf Proof}:  
To prove necessity, assuming that
$d(k) \le d_{\max}$  for all~$k$ with probability one, we must show the required condition on~$I(k)$ holds along all sample pathes. Indeed, letting~$k$ be given, since both \mbox{$d(k) \le d_{\max}$} and $d(k+1) \le d_{\max}$ with probability one, we claim this forces the required inequalities on $I(k)$. Without loss of generality, we provide a proof of the rightmost inequality for the case~{$I(k) \ge 0$} and note that a nearly identical proof is used for~{$I(k) < 0$}.  Indeed,  using the fact that~$X_{\min}$ is in the support~$\cal X$,  
there exists a suitably small neighborhood of~$X_{\min}$, call it $\mathcal{N}(X_{\min})$, such that
\[
P \left (X(k) \in \mathcal{N}(X_{\min}) \right) >0.
\]
Thus, given any arbitrarily small~$\varepsilon >0$, there exists some point $X_{\varepsilon}(k) < 0$ such that $
 X_{\varepsilon}(k) \in \mathcal{N}_\varepsilon(X_{\min}) ,
 $
  $|X_{\min} - X_\varepsilon(k)| < \varepsilon $ and leading to
 associated realizable loss $I(k)|X_
\varepsilon (k)|$.
 Noting that~$V_{\max}(k+1) = V_{\max}(k)$ and
\begin{align*}
 d(k + 1)= d(k) + \frac{{\left| {{X_\varepsilon (k) }} \right|I(k)}}{{{V_{\max }}(k)}} \leq d_{\max}
\end{align*}
we now substitute
\[{V_{\max }}(k) = \frac{{V\left( k \right)}}{{1 - d\left( k \right)}} >0
\]
into the inequality above and noting that $|X_\varepsilon(k)| \to |X_{\min}|$ as~$\varepsilon \to 0$, we obtain
\[
 I(k) \le \frac{{{d_{\max }} - d(k)}}{{\left( {1 - d\left( k \right)} \right)\left| {{X_{\min }}} \right|}}V(k).\]
To prove sufficiency,  we assume that the condition on  $I(k)$ holds along all sample pathes. We must show $d(k) \le d_{\max}$  for all~$k$ with probability one.  Proceeding by induction,  for~$k=0$, we trivially have~\mbox{$d(0)=0 \leq d_{\max}$} with probability one.  To complete the inductive argument, we assume that $d(k) \le d_{\max}$ with probability one, and  must  show~$d(k+1) \le d_{\max}$ with probability one. Without loss of generality, we again provide a proof for the case~{$I(k) \ge 0$} and note that a nearly identical proof is used for~{$I(k) < 0$}. Now, by noting that
\begin{align*}
d\left( {k + 1} \right) 
&= 1 - \frac{{V\left( {k + 1} \right)}}{{{V_{\max }}\left( {k + 1} \right)}},
\end{align*}
and $V_{\max}(k) \leq V_{\max}(k+1)$ with probability one, we split the argument into two cases: If $V_{\max}(k) < V_{\max}(k+1)$, then $V_{\max}(k+1) = V(k+1).$ Thus, we have $d(k+1)=0 \leq d_{\max}.$ On the other hand, if $V_{\max}(k) = V_{\max}(k+1)$, with the aid of the dynamics of account value, we have
\begin{align*}
d\left( {k + 1} \right) 
&= 1 - \frac{{V\left( k \right) + I\left( k \right)X\left( k \right)}}{{{V_{\max }}\left( k \right)}}\\
& \le 1 - \frac{{V\left( k \right) - I\left( k \right)\left| {{X_{\min }}} \right|}}{{{V_{\max }}\left( k \right)}}
\end{align*}
 Using the rightmost inequality condition on $I(k)$,  we obtain $ d\left( {k + 1} \right) \le {d_{\max }}$
which completes the proof. $\;\;\;\;\; \square$


\vspace{5mm}
{\bf Remarks}:
\label{Remarks of Drawdown Modulation Lemma}
First, it is worth to mentioning that since~\mbox{$d_{\max} < 1$}, any investment~$I(k)$ satisfying the inequality condition in the lemma assures survival. That is, along all sample pathes, $V(k) > 0$ for $k=0,1,2,...,N$. This is an easy consequence of the fact that at each stage $k$, drawdown of $100\%$ from some previous maximum for $V(k)$ never occurs.
Second, the satisfactions of the drawdown condition $d(k) \leq d_{\max}$ along sample pathes opens a door to solution of new drawdown-constrained optimization problems which involve parameters entering into an investment~$I(k)$ satisfying the inequality in the lemma; e.g., see Section~\ref{Section 4: Kelly-Based Optimization Problem} to follow. 
Finally, note that the i.i.d. assumption on the returns~$X(k)$ is not required in the proof of the lemma. What is needed is an explicit bound on the returns which is realizable with non-zero probability.

\vspace{5mm}
{\bf Feedback Control Realization}:
\label{Drawdown-Modulated Feedback Control}
 With the aid of the Drawdown Modulation Lemma, we can readily obtain a class of investment functions~$I(k)$ expressed as a linear feedback control parameterized by a gain $\gamma$ and leading to satisfaction of the drawdown specification. To be more specific, for each stage $k=0,1,2,...,N-1$, we define
\[
	M\left( k \right) \doteq \frac{{{d_{\max }} - d\left( k \right)}}{{1 - d\left( k \right)}}
\]
which we call {\it modulator function}. Now, using $M(k)$, we express~$I(k)$ in the feedback form
\[
	I(k) \doteq \gamma M(k) V(k)
\] 
with constraint 
$$
	-\frac{1}{ X_{\max}} \leq \gamma \leq   \frac{1} { |X_{\min}| }. 
$$ 
Note that the feedback gain $\gamma$ can be selected without regard for the modulator $M(k)$. This idea has a similar flavor to that of the celebrated Separation Theorem in linear control theory;~e.g., see~\cite{Chen_1995}. Henceforth, we call the investment $I(k)$ described above a {\it drawdown-modulated feedback controller} and the constraint on $\gamma$ above is denoted by writing $\gamma \in \Gamma.$

\vspace{5mm}

We note that this is not the only feedback-control realization which is possible. A more general class of feedback controls, to be pursued in our future work, can be formulated with a time varying gain~$\gamma(k)$. That is, one can equally well take~\mbox{$I(k) = \gamma(k)M(k)V(k)$} with $\gamma(k) \in \Gamma$ and still satisfy the drawdown requirement; see the conclusion of this paper for further discussion of promising research along these lines.

\vspace{8mm}
\section{Drawdown-Modulated Kelly Optimization}
\label{Section 4: Kelly-Based Optimization Problem}
\vspace{-2mm}
In this section, we consider one of many possible ways to incorporate the Drawdown Modulation Lemma into an optimization problem. To this end, we formulate a drawdown-constrained version of the so-called Kelly Optimization Problem.   To be more specific, for $k=0,1,2,...,N-1$,  assuming that the returns $X(k)$ are i.i.d. random variables with probability density function (PDF) denoted by $f_X$, we work with the drawdown-modulated feedback controller
$$
I(k) = \gamma M(k) V(k)
$$
with feedback gain $\gamma \in \Gamma$ treated as an optimization parameter.
 Now letting $V_\gamma(k)$ be the account value at stage $k$ induced by feedback gain $\gamma$,  the associated recursion is described by
\[
V_\gamma(k+1) = (1 + \gamma M(k) X(k) ) V_\gamma(k)
\]
and we seek to find $\gamma$ which achieves
$$
{g ^*} \doteq \mathop {\max }\limits_{\gamma  \in \Gamma } \frac{1}{N}\mathbb{E}\left[ {\log \frac{{{V_\gamma }\left( N \right)}}{{V\left( 0 \right)}}} \right].
$$
 Moreover, in the sequel, a corresponding maximizer is then denoted by
$
\gamma^*.
$

\vspace{5mm}
{\bf Remarks:}  
The limiting case obtained by letting~$d_{\max} \to 1$ leads to $M(k) \to 1$ which brings us back into the world of classical Kelly optimization;~e.g.,~see~\cite{Thorp_2006} and~\cite{Kelly_1956}. Thus, the classical Kelly problem can be viewed as a special case of our theory. In the next section, we consider a numerical example, based on historical data, to compare the trading performance obtained via drawdown-modulated feedback control with that obtained via a classical Kelly solution.  Given that this solution is often deemed to be too risky, some existing papers introduce the so-called {\it Fractional Kelly strategy} which was discussed in the introduction. One disadvantage of a fractional Kelly strategy is that it is often designed in an ad-hoc manner.  In contrast, the drawdown-modulated feedback control provides a systematic way to obtain a {\it time-varying}  problem solution.

\vspace{8mm}
\section{Numerical Example}
\label{Section 5: Illustrations}
\vspace{-2mm}
In this section, we provide a numerical example involving real stock data to illustrate how the drawdown modulation is used in practice. We consider the Kelly Optimization Problem described in the last section and use historical price data for Tesla Motors (ticker TSLA) covering the period December~31,~2013 to March~28,~2014; see Figure~\ref{fig:Stock_Price} where these underlying stock prices are plotted. The figure also shows an additional sixty-one days within the period from March~28,~2014 to June~24,~2014 of stock prices which will be used in an out-of-sample test described below. 

 \vspace{3mm}

\begin{center}
	\graphicspath{{Figs/}}
	\includegraphics[scale=0.42]{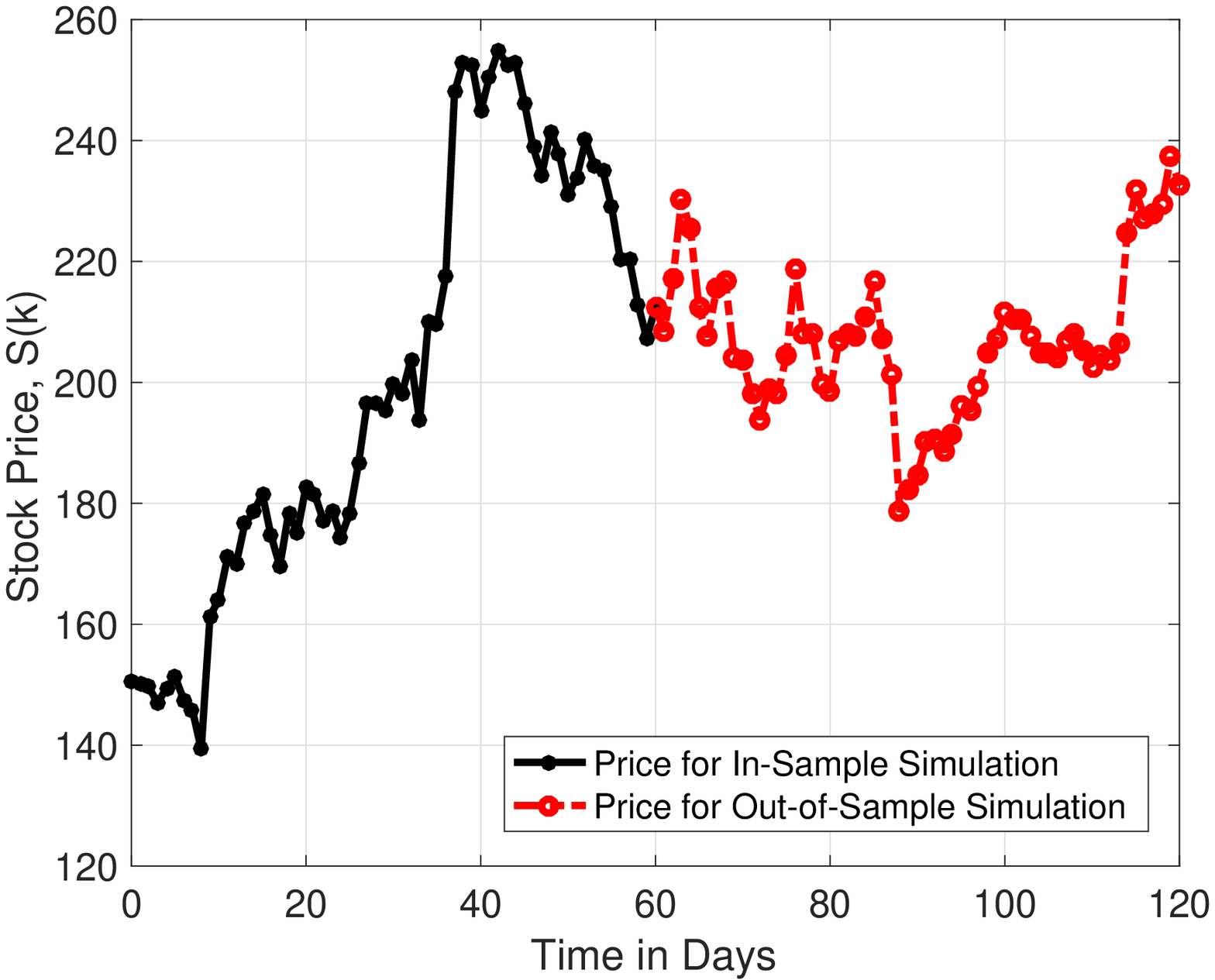}
	\figcaption{TSLA Stock Price, $S(k)$}
	\label{fig:Stock_Price}
\end{center}

\vspace{5mm}
Our goal is to first use the price data to estimate a probability mass function (PMF) model then calculate the optimal feedback gain $\gamma^*$ which maximizes the objective function 
$$
J(\gamma) \doteq \frac{1}{N} \mathbb{E}\left[\log \frac{V_\gamma(N)}{V(0)}\right]
$$  
with $V_\gamma(N)$ generated using the cash-financing constraint imposed.
We use this $\gamma^*$ to study the out-of-sample trading performance for the resulting controller. In this example we use $$
d_{\max} = 0.05
$$
as the maximum acceptable drawdown level. To be more specific, using the stock prices, we first calculate the corresponding returns
$$
X(k) = \frac{S(k+1)-S(k)}{S(k)}.
$$
 Now letting $x_i = X(i-1)$ denote the $i$-th calculated return for~$i=1,2,...,60$, we obtain the estimated PMF of the returns as the sum of impulses
\[
\hat{f}_X(x) = \frac{1}{60} \sum_{i=1}^{60} \delta (x - x_i)
\]
which is used as input to the optimization to be carried out. 
This PMF, plotted in Figure~\ref{fig:Stock_PMF}, has~$X_{\min} \approx -0.049$ and~\mbox{$X_{\max} \approx 0.157$.} Hence, the constraint set $\Gamma$ is described by 
$$
-6.354 \leq \gamma \leq 20.248
$$
for the Kelly Optimization Problem to be solved. 

\vspace{5mm}

Next, to evaluate~$J (\gamma)$ for each fixed $\gamma$, we perform Monte Carlo simulation to generate $100,000$  sample pathes for $S(k)$ from the PMF. Then we use these to estimate $\mathbb{E}[\log {V_\gamma }(N)]$ which is needed for the $J$ function evaluation. We note that use of the modulator automatically assures that the drawdown requirement is satisfied. Using the plot of the expected logarithmic growth in Figure~\ref{fig:Stock_Dmod}, we obtain a maximizer for $J(\gamma)$ given~by 
$$
\gamma^* \approx 11.15.
$$
\vspace{2mm}
 \begin{center}
 	\graphicspath{{Figs/}}
 	\includegraphics[scale=0.42]{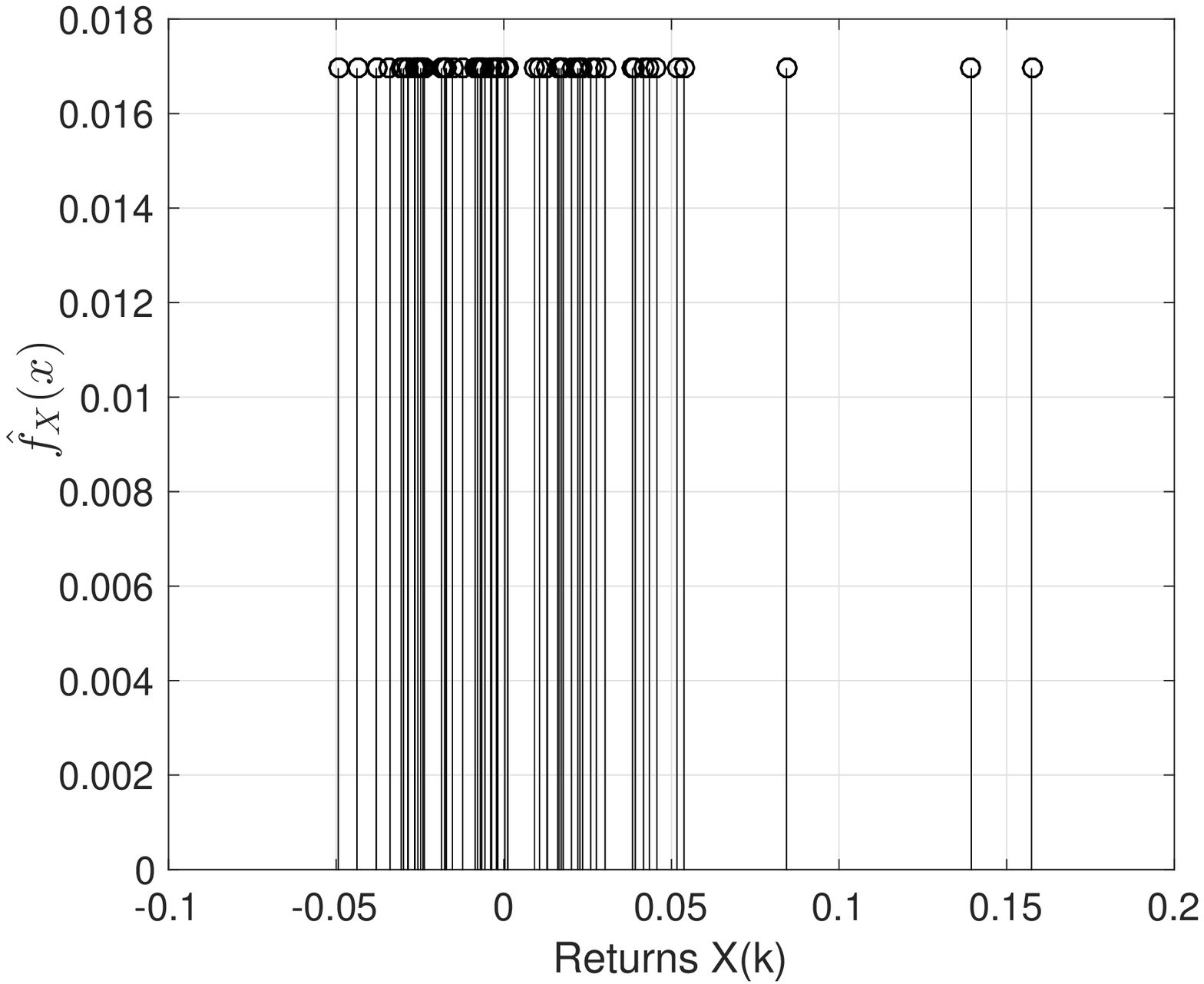}
 	\figcaption{Estimated PMF $\hat{f}_X(x)$ of Returns}
 	\label{fig:Stock_PMF}
 \end{center}
 
\vspace{5mm}
 
 \begin{center}
 	\graphicspath{{Figs/}}
 	\includegraphics[scale=0.42]{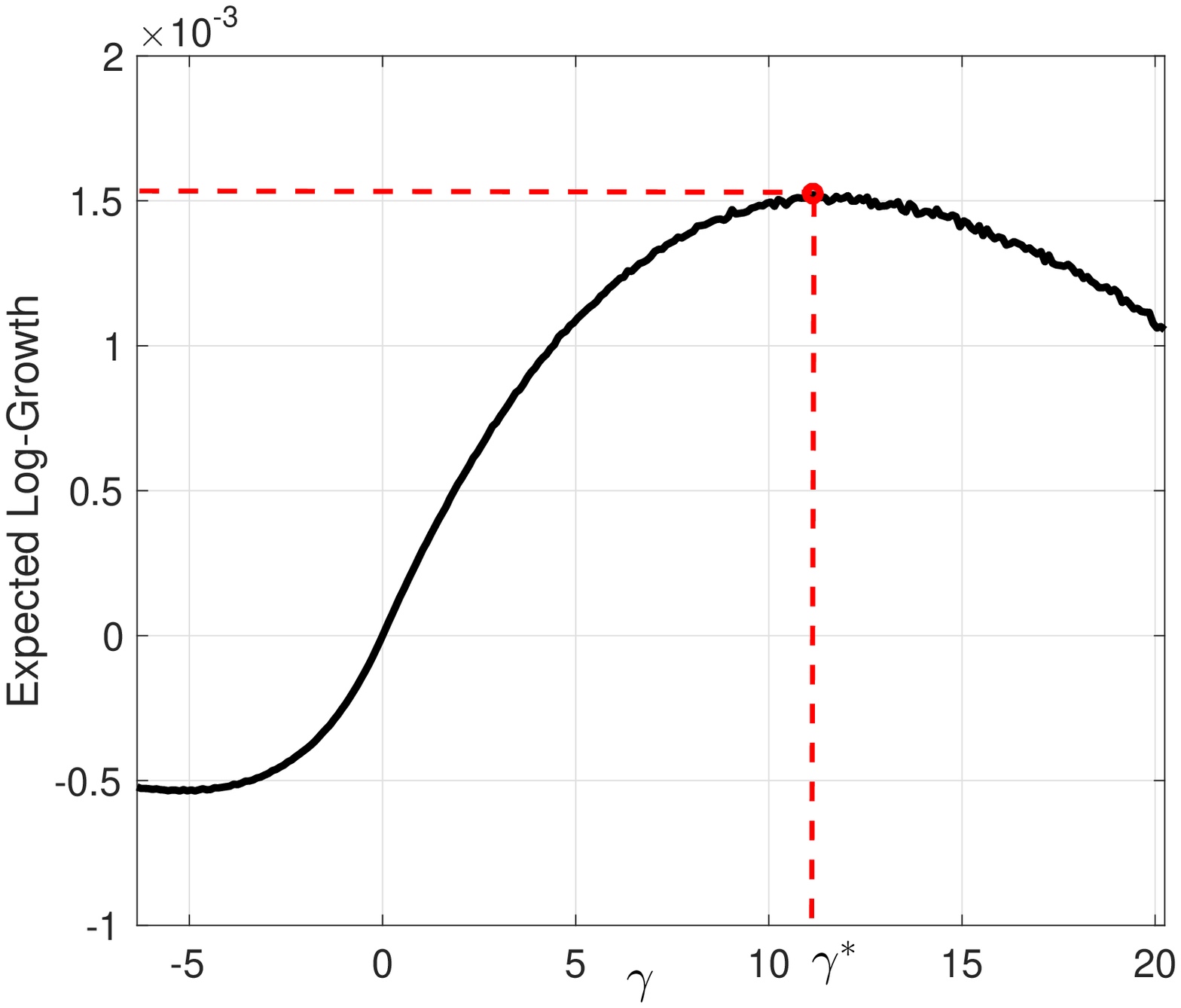}
 	\figcaption{Expected Growth Versus $\gamma$ via Drawdown Modulation}
 	\label{fig:Stock_Dmod}
 \end{center}

\vspace{5mm}
 {\bf In-Sample Trading Performance}: 
 \label{Section 5.1: In-Sample Simulation}
 Now, to examine the in-sample trading performance, we take the initial account value $V(0)=10,000$ and then compare the trading performance obtained via drawdown-modulated feedback control with that obtained via the classical Kelly solution. 
 Figure~\ref{fig:Controlled_V} depicts the sixty-one days in-sample simulation result. In particular, we see that the drawdown-modulated feedback control with $$
 \gamma^* \approx 11.15$$
  leads to  account value given by~$
 V(60)=1.146 \times 10^4 
 $
 with the overall percentage drawdown
 $$
 d_{\max}^* \approx 0.047
 $$
which is within the allowed upper limit of $5\%$. For comparison purposes, we also computed the classical Kelly solution which is obtained as a special case of our formulation as $d_{\max} \to 1$. In this case, the optimum turned out to be $\gamma^* = 1$ with the account value~\mbox{$
 V(60) = 1.412 \times 10^4
 $}
 and the overall percentage drawdown
 $
 d_{\max}^* \approx 0.186
 $
  which is much higher than the drawdown we obtained using modulation. 


\vspace{5mm}
 {\bf Out-of-Sample Trading Analysis}: 
 \label{Section 5.2: Out-of-Sample Simulation}
Next, we evaluate the out-of-sample performance for the second segment of stock prices given in Figure~\ref{fig:Stock_Price}. That is, we use the optimal feedback gain $\gamma^* \approx 11.15$ obtained previously to carry out sixty new trades.  Again, we begin with the same initial account value;~i.e., {$V(60) = 10,000$}. 
The associated trading performance is depicted in Figure~\ref{fig:Controlled_V_OutOfSample}. We see that the  the drawdown-modulated control leads to the terminal account value~\mbox{$
V(120) \approx    1.005\times 10^4
$}
with the overall percentage drawdown
$
d_{\max}^* \approx 0.05
$
as required.
In contrast, the classical Kelly strategy without a drawdown constraint leads to the terminal account value $
V(120) \approx 1.136 \times 10^4
$
and the overall percentage drawdown
$
d_{\max}^* \approx 0.225
$
which corresponds to $22\%$.
While this level is much higher than the specification~$d_{\max} = 0.05$ , the terminal account value is considerably higher. This higher return is to be expected since the drawdown risk was neglected. There is one interesting observation can be made from Figure~\ref{fig:Controlled_V_OutOfSample}. That is, when the maximum drawdown acceptable level $d_{\max}$ is hit, then the investment~$I(k)$ becomes zero. Thus, we see that the drawdown-modulated feedback control yields a ``stop-loss" type of behavior.

 \vspace{2mm}
 
 \begin{center}
 	\graphicspath{{Figs/}}
 	\includegraphics[scale=0.42]{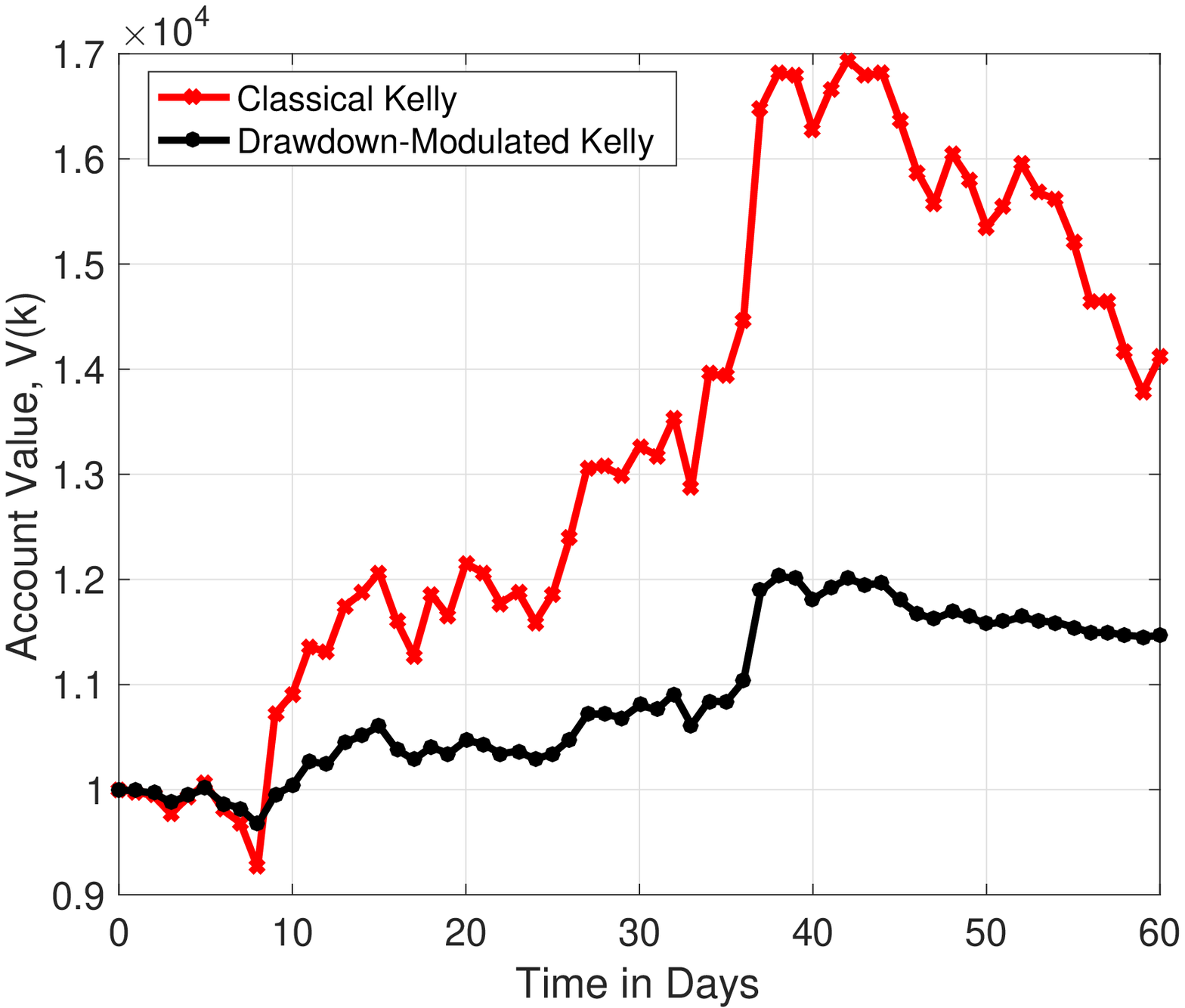}
 	\figcaption{Account Value Under Two Strategies (In-Sample)}
 	\label{fig:Controlled_V}
 \end{center}
 
\vspace{5mm}
 
\begin{center}
	\graphicspath{{Figs/}}
	\includegraphics[scale=0.42]{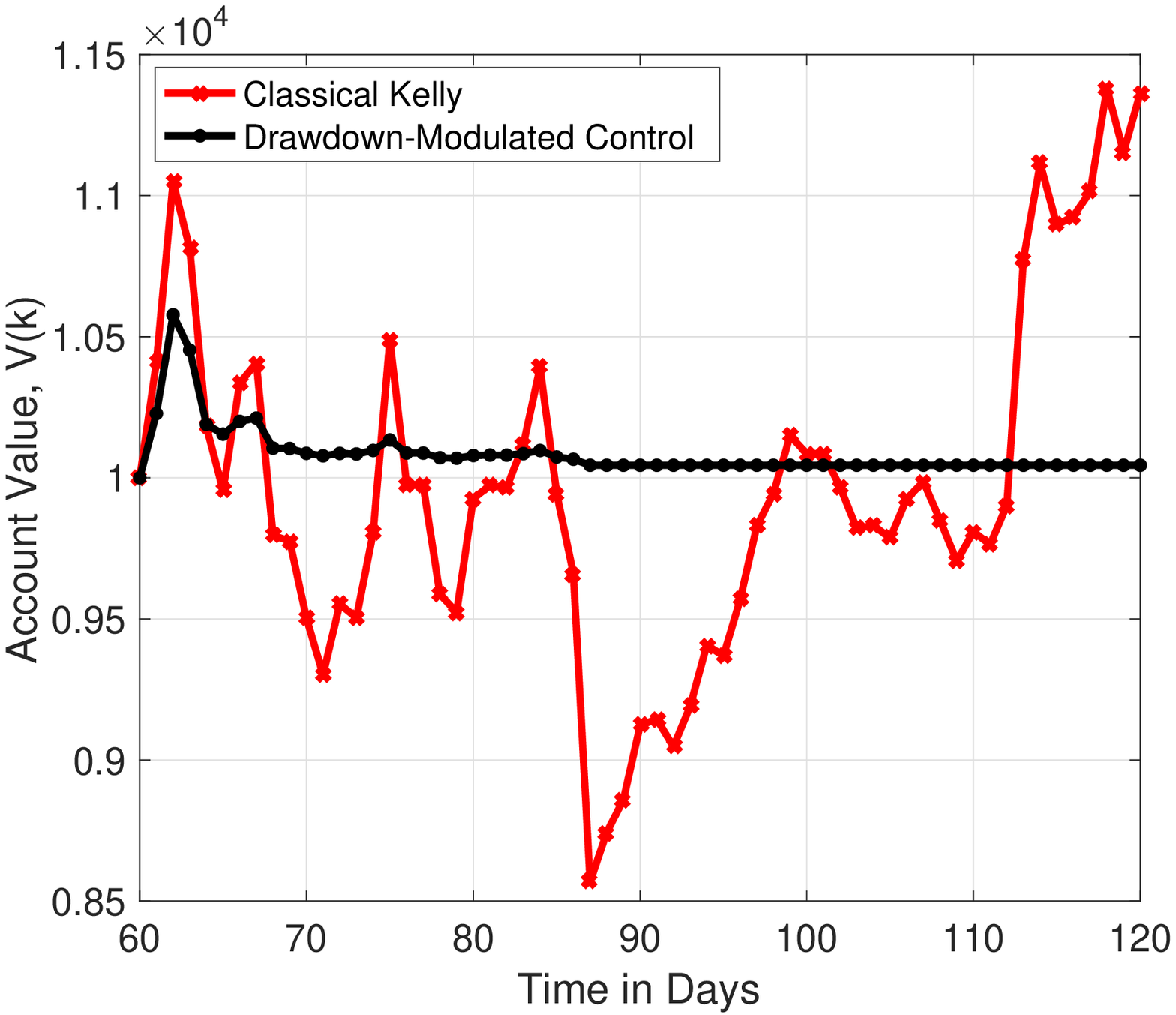}
	\figcaption{Account Value Under Two Strategies (Out-of-Sample)}
	\label{fig:Controlled_V_OutOfSample}
\end{center}

\vspace{8mm}
\section{Generalized Drawdown Modulation}
 \label{Section 6: Embellishment of  Drawdown Modulation Lemma: Portfolio Case }
\vspace{-2mm}
The Drawdown Modulation Lemma in Section~\ref{Section 3: The Drawdown Modulation Theorem} can be generalized to address a portfolio of $n$ stocks as follows: With $S_i(k)$ being the $i$-th stock price, we form the returns
\[
X_i(k) \doteq \frac{S_i(k+1) - S_i(k)}{S_i(k)}
\]
for $i=1,2,...,n$
and vectorize it as 
$$
X(k) \doteq [ X_1(k) \; X_2(k)\; \cdots X_n(k) ]^T.
$$
As in the scalar case, we assume known bounds
$$
X_{\min,i} \le X_i(k) \le X_{\max,i}
$$
with
$$
-1< X_{\min,i} < 0 < X_{\max,i}
$$ 
for~$i=1,2,...,n$ and~$k=0,1,2,...,N-1$. 
We let 
$\mathcal{V}$ denote the~$2^n$ vertices of the hypercube defined by $X_{\min,i}$ and~$X_{\max,i}$ above and assume that every point $v \in \mathcal{V}$ is in the support~$\cal{X}$ of~$X(k)$. Henceforth, for simplicity of notation, we take $X_{\min}$ and $X_{\max}$ to be vectors with~$i$-th components~$X_{\min,i}$ and~$X_{\max,i}$, respectively
and we let~$|X_{\min}|$  denote the vector with  $i$-th component by~$|X_{\min,i}|$.
Now letting~$I_i(k)$ be the~\mbox{$i$-th} component of the investment vector $I(k)$,
the associated account value is updated~as
\[
V(k+1) = V(k) +  I^T(k) X(k).
\]
Now denoting the positive and negative parts of $I(k)$ component-wise by $I_i^+(k ) \doteq \max\{I_i(k),0\}$ and~\mbox{$I_i^-(k) \doteq \min \{I_i(k),0\}$} respectively, we are now prepared to provide a result for the general case of an $n$-stock portfolio. The proof which proceeds along similar lines to that of the lemma in Section~\ref{Section 3: The Drawdown Modulation Theorem} is included for the sake of completeness.


\vspace{5mm}
{\bf Generalized Drawdown Modulation Lemma}: {\it An investment function $I(\cdot)$ guarantees that the maximum acceptable drawdown  level~$d_{\max}$ or less with probability one if and only if for all $k$, the condition
\[\left| {X_{\min }^T} \right|{I^ + }(k) - X_{\max }^T{I^ - }(k) \le M(k)V(k)\]
is satisfied along all sample pathes.
}

\vspace{5mm}
{\bf Proof}: To prove necessity, assuming that
$d(k) \le d_{\max}$  for all~$k$ with probability one, we must show the required condition on~$I(k)$ holds along all sample pathes. Indeed, letting~$k$ be given, since both \mbox{$d(k) \le d_{\max}$} and $d(k+1) \le d_{\max}$ with probability one, we claim this forces the required inequalities on~$I(k)$. Indeed, with the $i$-th component $I_i(k)$ being either~$I_i^+(k)$ or~$I_i^-(k)$, there exists a vertex~$v \in \mathcal{V}$   such that
	$$
	 v^T I(k)  = X^T_{\max} I^-(k) + X^T_{\min} I^+(k) \leq 0.
	$$ 
Using the fact that $v$ is  in the support~$\cal X$, it follows that there exists a neighborhood of $v$, call it $\mathcal{N}(v)$, such that
\[
P(X(k) \in \mathcal{N}(v)) > 0.
\]
Hence, given any arbitrarily small~$\varepsilon >0$, there exists some point $X_{\varepsilon}(k)$ such that 
$
X_{\varepsilon}(k) \in \mathcal{N}_\varepsilon(v)$,  $ ||v - X_\varepsilon(k)|| < \varepsilon $ and leading 
 to the  realizable loss
$
 X^T_\varepsilon (k) I(k) \leq 0.
$
Noting that we have
$
V_{\max}(k+1) = V_{\max}(k).
$
Thus, it follows that
\begin{align*}
 d(k + 1) = d(k) - \frac{ X^T_\varepsilon (k) I(k) }{{{V_{\max }}(k)}} \leq d_{\max }.
\end{align*}
Now, substituting 
\[
{V_{\max }}(k) = \frac{{V\left( k \right)}}{{1 - d\left( k \right)}} >0
\]
into the inequality above and noting that~\mbox{$X^T_\varepsilon (k) I(k)  \to v^TI(k)$} as~$\varepsilon \to 0$, we obtain
\[
|X_{\min}^T| I^+(k) - X_{\max}^T I^-(k)  \leq M(k)V(k).
\]
To prove sufficiency, we assume the condition on $I(k)$ holds along all sample pathes. We must show $d(k) \le d_{\max}$  for all~$k$ with probability one. Proceeding by induction,  for~$k=0$, we trivially have~$d(0)=0 \leq d_{\max}$ with probability one.  To complete the inductive argument, we assume that $d(k) \le d_{\max}$ with probability one, and must  show $d(k+1) \le d_{\max}$ with probability one. Now,  by noting that 
\begin{align*}
d\left( {k + 1} \right) 
&= 1 - \frac{{V\left( {k + 1} \right)}}{{{V_{\max }}\left( {k + 1} \right)}}
\end{align*}
and $V_{\max}(k) \le V_{\max}(k+1)$, we split the argument into two cases: If $V_{\max}(k) < V_{\max}(k+1)$, then $V_{\max}(k+1) = V(k+1).$ Thus, we have $d(k+1)=0 \leq d_{\max}.$ On the other hand, if $V_{\max}(k) = V_{\max}(k+1)$, with the aid of the dynamics of account value, we have
\begin{align*}
d\left( {k + 1} \right) 
&= 1 - \frac{{V\left( k \right) + I^T\left( k \right)X\left( k \right)}}{{{V_{\max }}\left( k \right)}}\\
& \le 1 - \frac{{V\left( k \right) - \left( {\left| {X_{\min }^T} \right|{I^ + }(k) - X_{\max }^T{I^ - }(k)} \right)}}{{{V_{\max }}\left( k \right)}}.
\end{align*}
 Using the given inequality condition on $I(k)$,  we obtain
$ d\left( {k + 1} \right) \le {d_{\max }}$
which completes the proof. $\;\;\;\;\; \square$

%

\vspace{5mm}
{\bf Remarks}: The drawdown-modulated feedback realization described for the single-stock case is generalized to this multi-stock case as follows: For the $i$-th stock, we use feedback gain~$\gamma_i$ and take the investment to be 
\[
I_i(k) \doteq \gamma_i M(k) V(k).
\]
Then, the condition in the Generalized Drawdown Modulation Lemma above leads to the constraint on $\gamma$ as follows:
\[
| X_{\min }^T| \gamma^+ -X_{\max }^T \gamma^- \leq 1
\]
where $\gamma^-$ and $\gamma^+$ have $i$-th component\[\gamma _i^ -  \buildrel\textstyle.\over= \min \{ {\gamma _i},0\} ;\;\;\; \gamma _i^ +  \buildrel\textstyle.\over= \max \{ {\gamma _i},0\} .\]
%

\vspace{8mm}
\section{Conclusion and Future Research}
\label{Section 7: Conclusion}
\vspace{-2mm}
In this paper we  introduced a new scheme, which we called drawdown-modulated feedback control. It enables us to  express the investment $I(k)$ as a linear feedback realization with a gain~$\gamma$ which leads to satisfaction of a given percentage drawdown specification with probability one. We also provided an illustration of our theory in the context of the Kelly Optimization Problem. The resulting drawdown-modulated feedback control which we obtain provides a systematic way to obtain a time-varying fractional Kelly strategy which takes care of the drawdown requirements. 

\vspace{5mm}
To further pursue this research, one obvious problem to consider would be the portfolio optimization version of the expected logarithmic growth maximization described in Section~\ref{Section 4: Kelly-Based Optimization Problem} with consideration of drawdown along the lines of Section~\ref{Section 6: Embellishment of  Drawdown Modulation Lemma: Portfolio Case }.  To find an associated optimal solution vector $\gamma^*$ in the portfolio scenario,  some efficient algorithm aimed at dealing with potential computational complexity issue might need to be used. 
 
\vspace{5mm}
 
Further to the optimization aspects of the work in this paper, we draw attention to the fact that the solution~$\gamma^*$ we obtained  is simply a pure gain. However, it can be argued that such a pure gain~$\gamma$ is not necessarily the true optimum for many problems. Recalling the discussion in Section~\ref{Drawdown-Modulated Feedback Control}, it may prove to be the case that a time-varying feedback gain $\gamma(k)$ may lead to superior performance.

\vspace{8mm}

                            

\end{document}